\documentclass[smallextended,envcountsect,]{svjour3}
\smartqed
\usepackage{graphicx}

\usepackage{latexsym, amsmath,amsxtra, amssymb, latexsym, amscd}
\usepackage{enumerate,cite,color}
\usepackage{hyperref}

\def\Limsup{\mathop{{\rm Lim}\,{\rm sup}}}

\def\R{\mathbb{R}}
\def\N{\mathbb{N}}
\def\Sol{\mbox{\rm Sol}}

\begin{document}

\title{On the existence  and the  stability of solutions in   nonconvex vector optimization$^\dagger$}

\titlerunning{On the existence  and the  stability of solutions in nonconvex   vector optimization}

\authorrunning{T.V. Nghi, L.N. Kien, N.V. Tuyen}


\author{Tran Van Nghi \and  Le Ngoc Kien  \and Nguyen Van Tuyen}

\institute{Tran Van Nghi   \at
	Department of Mathematics, Hanoi Pedagogical University 2,\\
	Xuan Hoa, Phu Tho, Vietnam \\
	tranvannghi@hpu2.edu.vn
           \and Le Ngoc Kien \at
           Faculty of Fundamental Science, Vietnam-Hungary Industrial University,
           \\
           Tung Thien, Hanoi, Vietnam;
           \\
           Department of Mathematics, Hanoi Pedagogical University 2,\\
           Xuan Hoa, Phu Tho, Vietnam \\
           lengockien@viu.edu.vn
           \and
           Nguyen Van Tuyen,  Corresponding author \at
           Department of Mathematics, Hanoi Pedagogical University 2,\\
           Xuan Hoa, Phu Tho, Vietnam \\
           nguyenvantuyen83@hpu2.edu.vn; tuyensp2@yahoo.com              
}

\date{Received: date / Accepted: date}

\maketitle
{\centerline{\em$^\dagger$Dedicated to Professor Q.H. Ansari on the occasion of his 65th birthday}}
\vskip 0.5cm
\begin{abstract}
The paper is devoted to the existence of weak Pareto solutions and  the weak sharp minima at infinity property for a general class of  constrained nonconvex vector optimization problems with unbounded constraint set via asymptotic cones and generalized asymptotic functions. Then we  show that these conditions are  useful for studying the solution stability of nonconvex vector optimization problems with linear perturbation.  We also provide some  applications for a subclass of  robustly quasiconvex vector optimization problems. 
\end{abstract}
\subclass{90C29 \and 49J30 \and  90C31 \and 90C26 \and  49J52
\\---------------------------------
\\
Communicated by Hong-Kun Xu.
}
\section{Introduction}

Vector optimization   is an important part of mathematical programming. It is used to denote a type of optimization problems, where two
or more objectives are to be minimized over a given set.
Vector optimization problems arise naturally in a variety of
theoretical and practical problems and have a wide range of applications. These problems  are important tools for problems in economics, management, multicriteria design optimization, water resource planning, medicine, mathematical biology, etc; see, for example, \cite{Ansari-2,Ansari-3,luc89,D-1959,Ehrgott,Khan-et al,Mor06a}.  
\medskip

An important and fundamental issue in theory of vector optimization problems is the investigation of existence and stability of  solutions; see, for example, \cite{Ansari-1,Ehrgott,luc89,Khan-et al,Borwein,Jahn,Luc-89-2,Chuong-24,Chuong-24-2,Huy-et-al-17}.  In this
paper, we consider the following vector optimization problem 
\begin{equation}\label{problem-1}
	\mathrm{Min}_{\R^m_+}\,\{f(x):=(f_1(x), \ldots, f_m(x))\;:\; x\in X\},\tag{VP}
\end{equation}
where  $X$  is a  closed and unbounded subset in $\R^n$ and $f_i\colon \R^n\to \R$, $i~=~1, \ldots, m$, are lower semicontinuous functions.
\medskip

Some authors have studied  theorems of Frank--Wolfe type  for polynomial vector optimization problems.
Kim et al. \cite{c7} showed the existence of Pareto  solutions of an unconstrained polynomial vector optimization problem under the assumptions that the
Palais--Smale type condition holds and the image of the objective function has a
bounded section. For the case where $X$ is a convex semi-algebraic set and $f$ is convex, Lee et al.
\cite{c3} proved that \eqref{problem-1} has a Pareto solution if and only if $f(X)$ has a
nonempty bounded section. Recently, Flores-Baz\'an et al. \cite{F-L-Vera}  introduced two notions of vector asymptotic functions for vector functions with
their image space ordered by the nonnegative orthant and  established coercivity properties, coercive and noncoercive existence results for solutions of
vector optimization problems. 
\medskip

The main difficulty in investigation the existence of (weak) Pareto solutions is due to directions in unbounded constraint set along which the objective function may decrease.  Asymptotic tools have proved to be very useful in order to
overcome the above difficulty; see, for example, \cite{Rele-Nedic-24,Stein,Ait,F-L-Vera,Lara-Tuyen-Nghi,AT}.   The
idea of the asymptotic  approach is based on the asymptotic
behavior of sets and functions.  When the
constraint set is unbounded, asymptotic tools
allow us to obtain coercivity properties and existence results for vector optimization problems.
\medskip

In this paper,  we investigate  the existence and the stability of weak Pareto solutions of constrained nonconvex vector optimization problems by using the tool of asymptotic analysis.  We make the following contributions to vector optimization problems:
\medskip

$\bullet$ As a first result, we present new results on existence of weak Pareto solutions and  the weak sharp minima at infinity property for a general class of  constrained nonconvex vector optimization problems with unbounded constraint set via asymptotic cones and generalized asymptotic functions.
\medskip

$\bullet$ We obtain stability results, including: the nonemptiness and compactness, the upper/lower semicontinuity,  and the closedness of the weak Pareto  solution map of the considered problem under linear perturbations.
\medskip

$\bullet$ We  introduce the concept of robustly quasiconvex for vector-valued functions
and  apply our previous results to the particular case when the 
objective function  is quasiconvex or $\alpha$-robustly
quasiconvex. By using the $q$-asymptotic function, we develop sufficient conditions for
the mentioned problems under  weaker assumptions.
\medskip

This paper is organized as follows. Section \ref{section2} provides the necessary background required for the subsequent development.  In Section \ref{section3}, we investigate the solution existence and the weak sharp minima at infinity.  Section \ref{section4} presents  stability results under a linear perturbation. Section \ref{section5} proposes  some applications  to the particular case when the 
objective function  is quasiconvex or $\alpha$-robustly
quasiconvex.
A summary of conclusions is explained in the last section.

\section{Preliminaries}\label{section2}

Throughout the paper, the considered  spaces are finite-dimensional Euclidean with the inner product and the norm being denoted, respectively, by $\langle \cdot, \cdot \rangle$ and  $\|\cdot\|$. The open 
ball centered at the origin with a   radius of $\alpha>0$ and the set of all positive integer numbers are  denoted, respectively, by   $\mathbb{B}_\alpha$ and  $\N$.  The nonnegative orthant of $\R^m$ is defined by $$\R^m_+:=\{x=(x_1, \ldots, x_m)\in\R^m\;:\; x_i\geq 0, i=1, \ldots, m\}.$$

Let  $\varphi\colon\mathbb{R}^{n}\rightarrow \overline{\mathbb{R}} :=
\mathbb{R}\cup \{ \pm \infty \}$ be an extended-real-valued  function.  The  {\em domain} and {\em epigraph} of $\varphi$ are defined, respectively, by
\begin{align*}
	\operatorname{dom}\varphi:=\{x\in \mathbb{R}^{n}\,:\,\varphi(x)<+\infty \} \ \ \text{and}
	\\
	\mathrm{epi}\,\varphi:=\{(x,t)\in \mathrm{dom}\,\varphi\times \mathbb{R}\,:\,\varphi(x)\leq
	t\}.
\end{align*}
We say that $\varphi$ is proper if $\operatorname{dom}\varphi\neq\emptyset$ and $\varphi(x)>-\infty$ for every $x\in \mathbb{R}^{n}$. We adopt the usual convention $\inf \nolimits_{\emptyset}\varphi:=+\infty$ and $\sup \nolimits_{\emptyset} \varphi:= -\infty$.

For a given $\alpha \in \mathbb{R}$,	we denote by 
$$\mathrm{lev}(\varphi, \alpha) := \{x\in \mathbb{R}^{n}:~\varphi(x)\leq \alpha \}$$ 
the sublevel set of $\varphi$  at value $\alpha$. For  a nonempty subset $X$ in $\mathbb{R}^n$, the set of all minimum points of $\varphi$ on $X$ is defined by
$$\mathrm{argmin}_{X}\varphi:=\{x\in X:~\varphi(x)\leq \varphi(y)\ \ \forall y\in X\}.$$
\begin{definition}\rm
	We say that the function $\varphi$ is:
	\begin{enumerate}[\rm(i)]
		\item {\em convex} if $\mathrm{epi}\,\varphi$ is a convex subset of $\R^n\times\R$.
		\item {\em quasiconvex} if $\mathrm{lev}(\varphi, \alpha)$ is a convex subset of $\R^n$ for all $\alpha\in\R$.  
	\end{enumerate}
\end{definition}
It is well-known that the function $\varphi$ is convex if and only if its domain is convex and for every $x, y \in {\rm dom}\,\varphi$
$$\varphi(\lambda x+(1-\lambda)y) \leq \lambda \varphi(x) + (1-\lambda) \varphi(y) \ \ \forall   \lambda \in [0, 1]$$ 
and it is quasiconvex if and only if 
for every $x, y \in {\rm dom}\,\varphi$,
$$\varphi(\lambda x + (1 - \lambda)y) \leq \max\{\varphi(x), \varphi(y)\}\ \ \forall  \lambda \in [0,1].$$	
Clearly, every convex function is also quasiconvex, but not vice versa. For example, the function $\varphi(x) := \sqrt{|x|}$ for all $x\in\R$ is quasiconvex but it is nonconvex.   

In \cite[Lemma 1.1]{PA}, the authors showed that the quasiconvexity is not stable under linear perturbations and hence, they introduced the so-called  $\alpha$-robustly quasiconvex functions.   

\begin{definition}[{\rm see \cite{BGJ,PA}}]\label{alpha:robust}   \rm
	For $\alpha \geq 0$, a proper function $\varphi: \mathbb{R}^{n} \rightarrow
	\overline{\mathbb{R}}$ is said to be {\it $\alpha$-robustly quasiconvex} if the 
	function $x \mapsto \varphi(x) + \langle u, x \rangle$ is quasiconvex for all $u \in
	\mathbb{B}_{\alpha}$. 
\end{definition}
As shown in \cite[p. 1091]{BGJ}, every convex function is also $\alpha$-robustly quasiconvex for all 	$\alpha \geq 0$, but not vice versa. Furthermore, every local minimum of an $\alpha$-robustly quasiconvex function  is a global one.  For a further study we refer the readers to \cite{BGJ,PA,PA-99}.

We now introduce the concept of robustly quasiconvex for vector-valued functions.
\begin{definition}\rm Let $f\colon\R^n\to\R^m$, $x\mapsto f(x):=(f_1(x), \ldots, f_m(x))$ be a mapping and $\alpha\geq 0$. We say that $f$ is:
	\begin{enumerate}[\rm(i)]
		\item {\em $\R^m_+$-quasiconvex} if for every $y=(y_1, \ldots, y_m)\in\R^m$, the  sublevel set 
		$$\mathrm{lev}(f, y):=\{x\in\R^n\,:\, f_i(x)\leq y_i, \ \ i=1, \ldots, m\}$$  
		of $f$ at $y$ is a convex subset of $\R^n$.  
		\item  {\em $\R^m_+$-$\alpha$-robustly quasiconvex} if the mapping $x\mapsto f(x)+\langle u, x\rangle$ is $\R^m_+$~-~quasiconvex for all $u=(u_1, \ldots, u_m)\in\R^{m\times n}$ with $\|u\|:=\max_{i=\overline{1, m}}\|u_i\|< \alpha$,  where
		$$\langle u, x\rangle:=(\langle u_1, x\rangle, \ldots, \langle u_m, x\rangle).$$ 
	\end{enumerate}  
\end{definition}
It is easy to see that $f$ is $\R^m_+$-quasiconvex (resp., $\R^m_+$-$\alpha$-robustly quasiconvex) if and only if each component $f_i$, $i=1, \ldots, m$, of $f$ is quasiconvex (resp., $\R^m_+$-$\alpha$-robustly quasiconvex).

Next, we recall some basic definitions and properties of asymptotic cones and functions, which can be found in \cite{AT}. 	
\begin{definition}\rm Let $X$ be a  subset of  $ \mathbb{R}^{n}$.
	\begin{enumerate}[\rm(i)] 
		\item The {\em asymptotic cone} of $X$, denoted by $X^\infty$, is the set defined by
		$$X^{\infty}:=\left \{  d\in \mathbb{R}^{n}~:~\exists~t_{k}\rightarrow +\infty, ~ \exists~x_{k}\in X,~\frac{x_{k}}{t_{k}}\rightarrow d\right \} \ \ \text{and}\ \ \emptyset^{\infty}=\emptyset.$$
		
		\item The (generalized) {\em recession cone} of $X$, denoted by $\mathrm{rec}\,X$, is the set defined by
		$$\mathrm{rec}\,X:=\{d\in \mathbb{R}^{n}~:~ x+td\in X\ \ \forall x\in X, \forall t\geq 0\} \ \ \text{and}\ \ \mathrm{rec}\,\emptyset=\emptyset.$$
	\end{enumerate}
\end{definition}
It follows from \cite[Propositions 2.1.1 and 2.1.2]{AT} that $X^\infty$ is a closed cone and $X$ is bounded if and only if $X^\infty=\{0\}$. While $\mathrm{rec}\,X$ is a convex cone, $\mathrm{rec}\,X\subset X^\infty$, and if $X$ is bounded, then $\mathrm{rec}\,X=\{0\}$.  When $X^\infty=\mathrm{rec}\,X$, then $X$ is called regular. In particular, if $X$ is convex, then it is regular and
\begin{equation*}\label{A1_convex} 
	X^{\infty}=\Big \{d\in \mathbb{R}^{n}:~x_{0}+t d\in
	X\  \forall~t \geq0\Big \} \, \text{ for any } x_{0}\in X,
\end{equation*}
see, for example, \cite[Proposition 2.1.5]{AT}.	
\begin{definition}\rm (see \cite{AT}) Let $\varphi\colon\mathbb{R}^n\to\overline{\mathbb{R}}$ be a proper function. The	{\em asymptotic function} $\varphi^{\infty}~:~\mathbb{R}^{n}\rightarrow ~\overline{\mathbb{R}}$ of   $\varphi$ is the function defined by
	\begin{equation*}\label{usual:formulas}
		\varphi^{\infty}(d):=\inf \left \{  \liminf_{k\rightarrow+\infty}\frac{\varphi(t_{k}d_{k})}{t_{k}}:~t_{k}\rightarrow+\infty,~d_{k}\rightarrow d\right \} \ \ \forall d\in\R^n.
	\end{equation*}
\end{definition}
By definition, it is easy to see  that $\varphi^\infty$ is lower semicontinuous and positively homogeneous, and the following relation holds $\operatorname{epi}\varphi^{\infty}=(\operatorname{epi}\varphi)^{\infty}.$ Moreover, when $\varphi$ is lower semicontinuous   and convex, then we have
\begin{equation*}
	\varphi^{\infty}(d)=\sup_{t>0}\frac{\varphi(x_{0}+td)-\varphi(x_{0})}{t}=\lim_{t\rightarrow
		+\infty}\frac{\varphi(x_{0}+td)-\varphi(x_{0})}{t}  
\end{equation*}
for all $x_{0} \in \operatorname{dom}\varphi$; see \cite[Proposition 2.5.2]{AT}.	 Recall that the function $\varphi$ is said to be \textit{lower semicontinuous} (lsc) at some $\bar x \in \R^n$  if for every real number $r < \varphi(\bar x)$ there exists a neighborhood $V$ of  $\bar x$ such that $r< \varphi(x)$ for all $x\in V$. 	We say that $\varphi$ is lsc around $\bar x$ if it is lsc at any point of some neighborhood of $\bar x$. If $\varphi$ is lsc at every $x \in X\subset \R^n,$ then $\varphi$ is said to be lsc on $X$.  By definition, it is clear that $\varphi$ is automatically lsc at $\bar x$ when $\varphi (\bar x) = -\infty$; when $\varphi (\bar x) = +\infty$ the lower semicontinuity of $\varphi$ means that the values of $\varphi$ remain as large as required, provided that one stays in some small neighborhood of $\bar x$.

It is well-known that, if $\mathcal{K}(\varphi)=\{0\}$, where $$\mathcal{K}(\varphi):=\{d\in\R^n\,:\, \varphi^\infty(d)\leq 0\},$$ 
then $\varphi$ is coercive, i.e.,
\begin{equation*}
	\lim_{\lVert x \rVert \to + \infty} \varphi(x) = + \infty.
\end{equation*}
Furthermore, if $\varphi$ is convex and lsc, then 
\begin{equation*}
	\mathcal{K}(\varphi)=\{0\}	\Longleftrightarrow	\varphi\mathrm{~is~coercive}\Longleftrightarrow \mathrm{argmin}_{\mathbb{R}^{n}}~\varphi \ \ \text{is nonempty and compact}
\end{equation*}
see \cite[Proposition 3.1.3]{AT}. However, when $\varphi$ is nonconvex, the asymptotic function $\varphi^{\infty}$ is not good enough 	for providing adequate information
on the behavior of $\varphi$ at infinity. For example,   let  $\varphi: \mathbb{R} \rightarrow \mathbb{R}$ be the function defined by $\varphi(x) = \sqrt{\lvert x \rvert}$. Then $\varphi$ is coercive and ${\rm argmin}_{\R}\,\varphi 
= \{0\}$, but $\mathcal{K}(\varphi) = \mathbb{R}$. Therefore, several authors  	have  proposed different notions for dealing with, specially, quasiconvex functions, see 
\cite{FFB-Vera,HLM,HLL,IL-1,Lara-4}. 
\begin{definition}\rm (see \cite{FFB-Vera,HLM}) 
	Let $\varphi\colon\mathbb{R}^n\to\overline{\mathbb{R}}$ be a proper function. The {\em $q$-asymptotic function}   of $\varphi$ is the function $\varphi^{\infty}_{q}: \mathbb{R}^{n} \rightarrow
	\overline{\mathbb{R}}$ given by
	\begin{equation*} 
		\varphi^{\infty}_{q} (d) := \sup_{x \in {\rm dom}\,\varphi} \sup_{t>0} \frac{\varphi(x+td) - \varphi(x)}{t} \ \ \forall d \in \mathbb{R}^{n}.
	\end{equation*}
\end{definition}	
As shown in \cite{FFB-Vera} that $\mathrm{epi}\,\varphi_q^\infty=\mathrm{rec}\,(\mathrm{epi}\,\varphi)$ and if   $\varphi$ is lsc and quasiconvex, then
\begin{equation*}\mathrm{argmin}_{\mathbb{R}^{n}}~\varphi\ \ \text{is nonempty and compact}~\mbox{ iff }
	~ \mathcal{K}_q(\varphi)=\{0\},
\end{equation*}
where 
$$\mathcal{K}_q(\varphi):=\{d\in\R^n\,:\, \varphi_q^\infty(d)\leq 0\}.$$ 
It can be shown that $\varphi_q^\infty$ is always positively homogeneous,  convex, and the following relation  holds  
\begin{equation*} 
	\varphi^{\infty}(d) \leq \varphi_{q}^{\infty}(d)\ \ \forall d\in\R^n,
\end{equation*}
see \cite{HLM}. Furthermore, the above inequality could be strict, see, for example,  \cite{Lara-4}. 

For a further study on generalized convexity and asymptotic analysis we refer the reader to \cite{Amara,att-but,Aus1,Aus2,AT,BGJ,Cambi,ffb1,ffb2,FFB-Vera,HKS,HLL,HLM,HL,IL-1,IL-3,IL-4,IL-5,Lara-1,Lara-4,lara-lopez,luc-pen,PA,Rele-Nedic-24,rock,rock-wet} and references therein.	

We now recall  definitions of the upper and lower semicontinuity to set-valued mappings. 	Let  $F : \mathbb{R}^n \rightrightarrows \mathbb{R}^m$ be a set-valued mapping and $\bar x\in\R^n$. The \textit{Painlev\'e--Kuratowski outer/upper limit} of $F$ at $\bar x$ is defined by
\begin{align*}
	\Limsup\limits_{x\rightarrow \bar x} F(x):=\bigg\{ y \in \mathbb{R}^m  :\; \exists x_k \rightarrow \bar x, y_k \rightarrow y \ \mbox{with}\ y_k\in F(x_k)\ \  \forall k=1,2,....\bigg\}.
\end{align*}

\begin{definition}\rm(see \cite{AT}) 
	We say that:
	\begin{enumerate}[\rm(i)]
		\item $F$ is {\em  upper semicontinuous} (usc henceforth) at $\bar{x}$ if, for any open 
		set $V \subset \mathbb{R}^m$ such that  $F(\bar{x}) \subset V$ there exists a neighborhood $U$ of $\bar{x}$ in $\R^n$  such that  $F(x) \subset V$ for all $x\in U$.
		\item  $F$ is {\em lower semicontinuous} (lsc) at $\bar x$ if $F(\bar x) \neq \emptyset$ and if, for any open set $V \subset \mathbb{R}^m$ such that $F(\bar x) \cap V \neq \emptyset$ there exists a neighborhood $U$ of $\bar x$  such that $F(x) \cap V \neq \emptyset$ for all $x \in U$.
		\item  $F$ is {\em continuous} at $\bar x$ if it is both usc and lsc at this point.
		\item $F$ is {\em closed} at $\bar x$  if for any sequences $x_k\in\R^n$ and $y_k\in\R^m$ one has
		\begin{equation*}
			x_k\to\bar x, y_k\to\bar y, y_k\in F(x_k) \Rightarrow \bar y\in F(\bar x).  
		\end{equation*}
		\item  $F$  is {\em locally bounded} at $\bar x$ if there exists a neighborhood $U$ of $\bar x$ such that $\cup_{x\in U} F(x)$ is bounded.
	\end{enumerate} 
\end{definition}	

\begin{remark}\rm 
	\begin{enumerate} [\rm(i)]
		\item By definition, it is easy to see that $F$ is closed at $\bar x$ if and only if the following condition holds
		\begin{equation*}
			\Limsup_{x \rightarrow \bar x} F(x)\subset F(\bar x).
		\end{equation*}
		\item By \cite[Theorem 1.4.1]{AT}, if $F$ is usc at $\bar x$ and $F(\bar x)$ is closed, then $F$ is closed at $\bar x$. Conversely, if $F$ is closed  and locally bounded   at $\bar x$,  then $F$ is usc at that point. We note here that the condition of locally boundedness is essential, i.e., if $F$ is closed at $\bar x$  but not locally bounded  at this point, then $F$ may not be usc at $\bar x$. For example, let $F\colon\R \rightrightarrows \R$ be defined by
		\begin{equation*}
			F(x)=
			\begin{cases}
				\{\tfrac{1}{x}\} \ \ &\text{if}\ \ x\neq 0,
				\\
				\{0\}  \ \ &\text{if}\ \ x= 0.
			\end{cases}
		\end{equation*} 
		Then, it is easy to see that $F$ is closed at $\bar x=0$ but $F$ is not usc at this point. 
	\end{enumerate}
	
\end{remark}

\section{The solution existence and the weak sharp minima at infinity}\label{section3}
In this    section, we consider the following vector optimization problem:
\begin{equation}\label{problem}
	\mathrm{Min}_{\R^m_+}\,\{f(x):=(f_1(x), \ldots, f_m(x))\;:\; x\in X\},\tag{VP}
\end{equation}
where  $X$  is a  closed and unbounded subset in $\R^n$, and $f_i\colon \R^n\to \R$, $i\in I$ with $I:=\{1, \ldots, m\}$, are lsc functions. The set
$$\Lambda:=\bigg\{\lambda\in \R^m_+\;:\; \sum_{i=1}^{m}\lambda_i=1\bigg\}$$
is called the unit simplex of $\R^m$. For each $\lambda\in \Lambda$, we consider the following scalarization problem
\begin{equation}\label{problem-lam}
	\min\,\bigg\{\varphi_\lambda(x):=\langle\lambda, f(x)\rangle\;:\; x\in X\bigg\}. \tag{VP$_\lambda$}
\end{equation}
The solution set of \eqref{problem-lam} is denoted by $S(\lambda)$.

\begin{definition}
	{\rm Let $\bar x\in X$. We say that:
		\begin{enumerate}[(i)]
			\item $\bar x$ is a {\em  weak Pareto solution} of \eqref{problem} if
			\begin{equation*}
				f(X)\cap(f(\bar x)-\mathrm{int}\,\R^m_+)=\emptyset.
			\end{equation*}
			The set of all weak Pareto solutions of \eqref{problem} is denoted by $\mathrm{Sol}^w\,\eqref{problem}$.
			\item  $\bar x$ is a {\em  Pareto solution} of \eqref{problem} if
			\begin{equation*}
				f(X)\cap(f(\bar x)-\R^m_+\setminus\{0\})=\emptyset.
			\end{equation*}
			The set of all  Pareto solutions of \eqref{problem} is denoted by $\mathrm{Sol}\,\eqref{problem}$.
		\end{enumerate}
	}
\end{definition}
\begin{remark}\rm  By definition, it is clear that
	\begin{equation*}
		\mathrm{Sol}\,\eqref{problem}\subset \mathrm{Sol}^w\,\eqref{problem}. 
	\end{equation*}
	
\end{remark}

The following result is  well-known but we give a proof here for the sake of the reader.
\begin{proposition}\label{Proposition 3.1} The following inclusion holds true
	\begin{equation*}
		\bigcup_{\lambda\in\Lambda} S(\lambda)\subset \mathrm{Sol}^w\eqref{problem}.
	\end{equation*}
\end{proposition}
\noindent{\it Proof} 
	Fix any $\lambda\in\Lambda$ and let $\bar x\in S(\lambda)$. We need to show that $\bar x\in \mathrm{Sol}^w\eqref{problem}$. If otherwise, then there exists $x_0\in X$ such that 
	$$f(x_0)\in f(\bar x)-\mathrm{int}\,\R^m_+,$$ 
	or, equivalently, $f_i(x_0)<f_i(\bar x)$ for all $i\in I$. This and the fact that $\lambda\in\Lambda$ imply that 
	$$\sum_{i\in I}\lambda_i f_i(x_0)<\sum_{i\in I}\lambda_i f_i(\bar x),$$
	contrary to $\bar x\in S(\lambda)$. The proof is complete. 
\qed 
\begin{theorem}[The coercivity and the weak sharp minima   at infinity]\label{weak-sharp-Thrm}  
	Assume that 
	\begin{equation}\label{equa-CQ}
		X^{\infty} \cap \left(\bigcup_{i\in I}\mathcal{K}(f_i)\right)=\{0\},
	\end{equation}
	where $\mathcal{K}(f_i):=\{d\in\mathbb{R}^n\,:\, f_i^\infty(d)\leq 0\}$, $i\in I$. 	Then the following statements hold:
	\begin{enumerate}[$(a)$]
		\item For each $i\in I$, $\mathrm{argmin}_X f_i(x)$ is nonempty and compact. Consequently, $f$ is bounded from below on $X$, i.e., there exists an element $a\in\R^m$ such that
		$$f(X)\subset a+\R^m_+.$$
		
		\item For each $\lambda\in\Lambda$, $S(\lambda)$ is nonempty and compact. Consequently,  			 $\mathrm{Sol}^w\,\eqref{problem}$ is nonempty.
		
		\item For all $\lambda\in\Lambda$, the function $\varphi_\lambda(\cdot):=\langle\lambda, f(\cdot)\rangle$ is coercive on $X$, i.e., \begin{equation*}
				\lim_{x \overset{X}{\longrightarrow} \, \infty} \varphi_\lambda(x)= + \infty,
		\end{equation*}
		where $x \xrightarrow{X} \infty$ means that $\|x\|\to \infty$ and $x\in X$.
		
		\item  $\mathrm{Sol}^w\,\eqref{problem}$ is compact.

		\item Problem \eqref{problem} has a weak sharp minimum at infinity, i.e., there exist constants $c>0$ and $R>0$ such that
		\begin{equation*} 
			\mathrm{dist}(f(x), f(\mathrm{Sol}^w\eqref{problem})) \geq c \, \mathrm{dist}(x, \mathrm{Sol}^w\eqref{problem}) \ \ \forall   x \in X \setminus \mathbb{B}_R. 
		\end{equation*}
	\end{enumerate}
\end{theorem}
\noindent{\it Proof} 
	$(a)$:  By definition, it is clear that $X^\infty$ and $\mathcal{K}(f_i)$, $i\in I$, are cones with apex at the origin. This and   \eqref{equa-CQ} imply that $X^\infty\cap \mathcal{K}(f_i)=\{0\}$ for all $i\in I$.  Hence, it follows from   \cite[Theorem 4.1]{Rele-Nedic-24} that $\text{argmin}_X f_i(x)$ is nonempty and compact for every $i\in I$. This implies that $f$ is bounded from below on $X$. 
	
	$(b)$: Fix any $\lambda\in\Lambda$ and let $\varphi_\lambda$ be the function defined by 
	$$\varphi_\lambda(x):=\sum_{i\in I}\lambda_i f_i(x) \ \ \forall x\in\R^n.$$
By \eqref{equa-CQ}, $f^\infty_i(d)>0$ for all $d\in X^\infty\setminus\{0\}$. This and \cite[Proposition 2.6.1(b)]{AT} imply that
	\begin{equation*}
		(\varphi_\lambda)^\infty(d)\geq \sum_{i\in I}\lambda_i(f_i)^\infty(d)>0 \ \ \forall d\in X^\infty\setminus\{0\}.
	\end{equation*}
	Consequently, $X^\infty\cap\mathcal{K}(\varphi_\lambda)=\{0\}$. Thus, for each $\lambda\in \Lambda$, $S(\lambda)$ is nonempty and compact  due to  \cite[Theorem 4.1]{Rele-Nedic-24}. This and Proposition \ref{Proposition 3.1} imply that $\mathrm{Sol}^w\,\eqref{problem}$ is nonempty.

	$(c)$: Fix any $\lambda\in\Lambda$. Then we have $X^\infty\cap\mathcal{K}(\varphi_\lambda)=\{0\}$.  By \cite[Theorem 3.1(c)]{Lara-Tuyen-Nghi},  $\varphi_\lambda$ is coercive on $X$.  
	
	$(d)$: To prove the compactness of $\mathrm{Sol}^w\,\eqref{problem}$, we need to show that $\mathrm{Sol}^w\,\eqref{problem}$ is closed and bounded. Firstly, we claim that $\mathrm{Sol}^w\,\eqref{problem}$  is bounded. Indeed, if otherwise, then there exists a sequence $x_k \in \mathrm{Sol}^w\,\eqref{problem}$ such that $x_k\to \infty$ as $k\to\infty$. By passing to subsequences if necessary we may assume that $d_k:=\frac{x_k}{\|x_k\|}$ converges to some $d\in\R^n$. Clearly, $d\in X^{\infty}$ and  $\|d\|=1$.  Fix any $x_0\in X$. Then for each $k\in\N$, due to the fact that $x_k \in \mathrm{Sol}^w\,\eqref{problem}$,   there exists $i_k\in I$ such that $f_{i_k}(x_k)\leq f_{i_k}(x_0)$. Since $i_k\in I$ for all $k\in N$, there  exists  $i_0\in I$ such that $i_k=i_0$ for infinite many $k$. Hence, by passing to a subsequence if necessary, we may assume that $f_{i_0}(x_k)\leq f_{i_0}(x_0)$
	for every $k\in \N$. Put $t_k:=\|x_k\|$. We have
	$$\frac{f_{i_0}(t_kd_k)}{t_k}\leq \frac{f_{i_0}(x_0)}{t_k}.$$
	This implies 
	$$\liminf_{k\rightarrow+\infty}\frac{f_{i_0}(t_{k}d_{k})}{t_{k}}\leq \liminf_{k\rightarrow+\infty} \frac{f_{i_0}(x_0)}{t_k} =0.$$
	It follows that $f_{i_0}^{\infty}(d)\leq 0$. This means that $d\in \mathcal{K}(f_{i_0})$. Hence 
	$$d\in X^{\infty} \cap \left(\bigcup_{i=1}^m\mathcal{K}(f_i)\right)\ \ \text{and}\ \ \|d\|=1,$$
	contrary to \eqref{equa-CQ}. Thus $\mathrm{Sol}^w\,\eqref{problem}$ is bounded. We now show that $\mathrm{Sol}^w\,\eqref{problem}$ is closed. Assume that $x_k$ is a sequence in  $\mathrm{Sol}^w\,\eqref{problem}$ such that $x_k$  converges to some $\bar x$. We claim that $\bar x\in \mathrm{Sol}^w\,\eqref{problem}$ and so $\mathrm{Sol}^w\,\eqref{problem}$ is closed. If otherwise, $\bar x\notin \mathrm{Sol}^w\,\eqref{problem}$. Hence, there exists $x_0\in X$ such that 
	\begin{equation}\label{equa-14}
		f_i(x_0)<f_i(\bar x) \ \ \ \forall i\in I.
	\end{equation}
	It follows from the lower semicontinuity of $f_i$  at $\bar x$ that
	\begin{equation*}
		f_i(\bar x)\leq \liminf_{k\to\infty}f_i(x_k) \ \ \ \forall i\in I.
	\end{equation*}
	This and \eqref{equa-14} give us that there exists $k\in\N$ large enough such that
	\begin{equation*}
		f_i(x_0)<f_i(x_k)  \ \ \forall  i\in I,
	\end{equation*}
	contrary  to the fact that $x_k\in \mathrm{Sol}^w\,\eqref{problem}$. Hence $\mathrm{Sol}^w\,\eqref{problem}$ is closed, as required.

$(e)$ Let $\psi\colon \mathbb{R}^n\to \mathbb{R}\cup\{+\infty\}$ be the achievement function of problem \eqref{problem} and defined by
\begin{equation*}
\psi(x):=\sup_{y\in X}\min_{i\in I}\left(f_i(x)-f_i(y)\right) \ \ \forall x\in\mathbb{R}^n.
\end{equation*}
Then, by \cite[Lemma 3.1]{Son-Ha-18} and the fact that $\mathrm{Sol}^w\eqref{problem}\neq\emptyset$, one has $\psi(x)<+\infty$ for all $x\in\mathbb{R}^n$ and 
\begin{equation*}
\mathrm{Sol}^w\eqref{problem}=\{x\in X\,:\, \psi(x)=0\}.
\end{equation*}

{\bf Claim 1.} {\em There exist $c_1>0$ and $R>0$ such that
\begin{equation}\label{equa-3-new}
\psi(x)\geq c_1\|x\| \ \ \forall x\in X\setminus\mathbb{B}_R.
\end{equation}  
}
  
Indeed, since $X^\infty\cap\mathcal{K}(f_i)=\{0\}$ for all $i\in I$, it follows from \cite[Theorem 7(c)]{Lara-Tuyen-Nghi} that $f_i$, $i\in I$, are coercive on $X$. By \cite[Proposition 3.1.2]{AT}, it is easy to see that there exist $\gamma>0$ and $R>0$ such that
\begin{equation*} 
f_i(x)\geq \gamma\|x\|\ \ \forall x\in  X \setminus \mathbb{B}_R.
\end{equation*}  
Then for any $x\in  X \setminus \mathbb{B}_R$ and $y\in X$, we have
\begin{equation*}
f_i(x)-f_i(y)\geq \gamma\|x\| - f_i(y) \ \ \forall i\in I.
\end{equation*}
Hence,
\begin{equation*}
\min_{i\in I}\left(f_i(x)-f_i(y)\right) \geq \gamma\|x\| +\min_{i\in I}(- f_i(y))
\end{equation*}
and we therefore get
\begin{equation*}
\sup_{y\in X}\min_{i\in I}\left(f_i(x)-f_i(y)\right) \geq \gamma\|x\|+\sup_{y\in X}\min_{i\in I}(- f_i(y)),
\end{equation*}
or, equivalently,
\begin{equation*}
\psi(x)\geq \gamma\|x\| -\inf_{y\in X}\max_{i\in I}f_i(y). 
\end{equation*}
By (a), there exists $a\in\mathbb{R}^m$ such that $f(X)\subset a+\mathbb{R}^m_+$. Hence,
\begin{equation*}
f_i(y)\geq \min_{k\in I} a_k 
\end{equation*}
for all $i\in I$ and $y\in X$. This implies that $\mathfrak{m}:=\inf_{y\in X}\max_{i\in I}f_i(y)\in\mathbb{R}$ and so   
\begin{align*}
\psi(x)&\geq \gamma\|x\| -\mathfrak{m}
\\
&=\frac{\gamma}{2}\|x\| +\left(\frac{\gamma}{2}\|x\|-\mathfrak{m}\right)\geq\frac{\gamma}{2}\|x\|  
\end{align*}  
for all $x\in X$ with $\|x\|$ large enough. Hence, by increasing $R$ (if necessary),  we may assume that 
\begin{equation*}
\psi(x)\geq \frac{\gamma}{2}\|x\| \ \ \forall x\in  X \setminus \mathbb{B}_R.
\end{equation*}
Then, the constant $c_1:=\frac{\gamma}{2}>0$ satisfies \eqref{equa-3-new}. 

\medskip
{\bf Claim 2.} {\em We have 
\begin{equation*}
\psi(x)\leq \mathrm{dist}(f(x), f(\mathrm{Sol}^w\eqref{problem})) \ \ \forall x\in\mathbb{R}^n.
\end{equation*}}
Indeed, fix any $\bar x\in \mathrm{Sol}^w\eqref{problem}$, then for every $x\in \mathbb{R}^n$ and $y\in X$, one has
\begin{align*}
f_i(x)-f_i(y)&=f_i(\bar x)-f_i(y) + f_i(x)-f_i(\bar x)
\\
&\leq f_i(\bar x)-f_i(y) +\max_{k\in I}\left(f_k(x)-f_k(\bar x)\right)
\\
&\leq f_i(\bar x)-f_i(y) +\|f(x)-f(\bar x)\| \ \ \forall i\in I.
\end{align*} 
Hence,
\begin{align*}
\min_{i\in I}\left(f_i(x)-f_i(y)\right) \leq \min_{i\in I}\left(f_i(\bar x)-f_i(y)\right)+\|f(x)-f(\bar x)\|.
\end{align*}
This implies that
\begin{equation*}
\sup_{y\in X}\min_{i\in I}\left(f_i(x)-f_i(y)\right)\leq\sup_{y\in X}\min_{i\in I}\left(f_i(\bar x)-f_i(y)\right) +\|f(x)-f(\bar x)\|,
\end{equation*}
or, equivalently,
\begin{align*}
\psi(x)\leq \psi(\bar x)+  \|f(x)-f(\bar x)\|.
\end{align*}
Since $\bar x\in \mathrm{Sol}^w\eqref{problem}$, $\psi(\bar x)=0$ and hence
\begin{equation*}
\psi(x)\leq  \|f(x)-f(\bar x)\|.
\end{equation*}
This implies that
\begin{equation*}
\psi(x)\leq \inf_{\bar x\in \mathrm{Sol}^w\eqref{problem}} \|f(x)-f(\bar x)\| =\mathrm{dist}(f(x), f(\mathrm{Sol}^w\eqref{problem})),
\end{equation*} 
as required.

\medskip
{\bf Claim 3.} {\em There exists  $R>0$ such that
\begin{equation*}
\|x\|\geq  \frac{1}{2}\,\mathrm{dist}\,(x, \mathrm{Sol}^w\eqref{problem}) \ \ \forall x\in\mathbb{R}^n, \|x\|\geq R.
\end{equation*}
}

Indeed, by the compactness of $\mathrm{Sol}^w\eqref{problem}$, we can choose $R$ large enough such that $\mathrm{Sol}^w\eqref{problem}\subset \mathbb{B}_R$. Take any $\bar x\in \mathrm{Sol}^w\eqref{problem}$ and $x\in\mathbb{R}^n$ with $\|x\|\geq R$, we have 
\begin{align*}
\|x\|&\geq \|x\|+\frac{1}{2}(\|\bar x\|-\|x\|)= \frac{1}{2}(\|x\|+\|\bar x\|)
\\
&\geq \frac{1}{2}\|x-\bar x\|\geq \frac{1}{2}\mathrm{dist}\,(x, \mathrm{Sol}^w\eqref{problem}),  
\end{align*}   
as required.

\medskip
We are now in a position to finish the proof of (e). Put $c:=\frac{c_1}{2}>0$.    Then, for $R$ large enough, we have
\begin{align*}
\mathrm{dist}(f(x), f(\mathrm{Sol}^w\eqref{problem}))&\geq \psi(x)\geq c_1\|x\|
\\
& \geq \frac{c_1}{2} \mathrm{dist}\,(x, \mathrm{Sol}^w\eqref{problem}) 
\\
&=c\, \mathrm{dist}\,(x, \mathrm{Sol}^w\eqref{problem})
\end{align*}
 for all $x\in X\setminus\mathbb{B}_R$ and which completes the proof.  
\qed

\begin{remark}\rm If condition \eqref{equa-CQ} is satisfied, then the Pareto solution set $\Sol\eqref{problem}$ is nonempty and bounded. Indeed, analysis similar to that in the proof of Proposition \ref{Proposition 3.1} shows that 
\begin{equation*}
	\bigcup_{\lambda\in\Lambda^\circ} S(\lambda)\subset \mathrm{Sol}\eqref{problem},
\end{equation*} 
where 
$$\Lambda^\circ:=\{\lambda\in\R^m: \lambda_i>0, i\in I, \sum_{i\in I}\lambda_i=1\}.$$
Combining this together with Theorem \ref{weak-sharp-Thrm}(b), we conclude that $\Sol\eqref{problem}$ is nonempty. The boundedness of $\Sol\eqref{problem}$ follows from Theorem \ref{weak-sharp-Thrm}(d) and the fact that $\Sol\eqref{problem}\subset \mathrm{Sol}^w\eqref{problem}$. 

However, the following example shows that condition \eqref{equa-CQ}  does not guarantee the closedness of  $\Sol\eqref{problem}$.  
\end{remark}
\begin{example}\label{Example 3.1}\rm Let $X\subset\R^2$ and $f\colon \R^2\to\R^2$ be defined, respectively, by
\begin{align*}
X:=[-1, 0]\times\{1\} &\cup \{x\in\R^2\,:\, x_1+x_2=1, 0\leq x_1\leq 1\}
\\
&\cup \{x\in\R^2\,:\, x_1-x_2=1, x_1\geq 1\}
\end{align*}
and $f(x)=(f_1(x), f_2(x)):=(x_1, x_2)$ for all $x=(x_1, x_2)\in\R^2$. Then, an easy computation shows that
$$\mathcal{K}(f_1)=\{d\in\R^2: d_1\leq 0\}, \mathcal{K}(f_2)=\{d\in\R^2: d_2\leq 0\},$$
and
$$X^\infty=\{d\in\R^2: d_1=d_2, d_1\geq 0\}.$$
Hence, $X^\infty\cap (\mathcal{K}(f_1)\cup\mathcal{K}(f_2))=\{0\}$, i.e., condition \eqref{equa-CQ} is satisfied.  However, it is easy to see that
$$\Sol\eqref{problem}=\{(-1, 1)\}\cup \{x\in\R^2: x_1+x_2=1, 0<x_1\leq 1\}$$
and this set  is not closed in $\mathbb{R}^2$.
\begin{center}
	\begin{figure}[htp]
		\begin{center}
			\includegraphics[height=7cm,width=9.5cm]{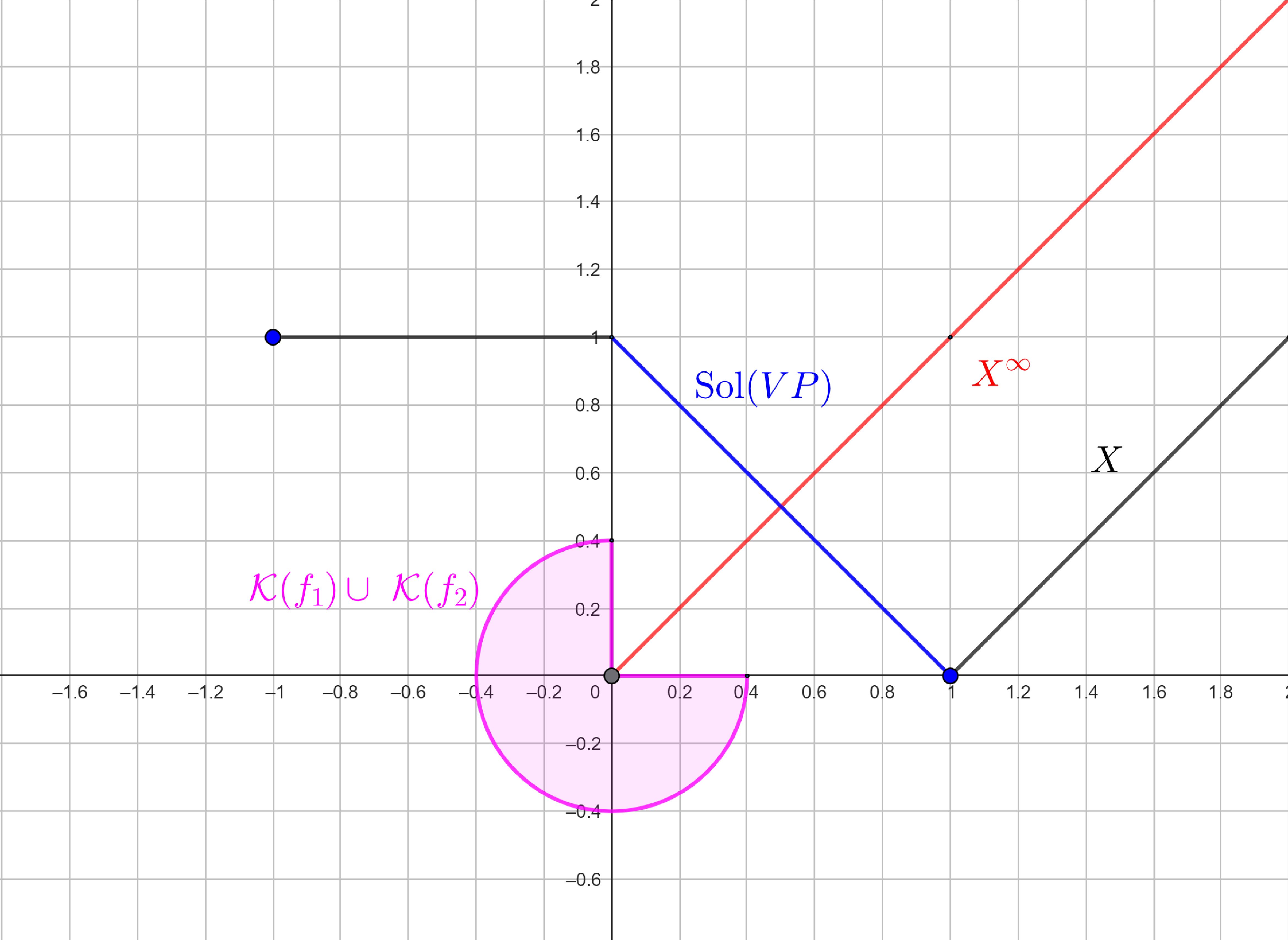}
		\end{center}
		\caption{Illustration of Example \ref{Example 3.1}.}\label{Fig_1}
	\end{figure}
\end{center}
\end{example}

\section{The solution stability results under linear perturbation}\label{section4}

For every $u=(u_1, \ldots, u_m)\in\R^{m\times n}$, we define the function $f^u\colon\R^n\to \R^m$ 
by $$f^u(x):=f(x)-\langle u, x\rangle:=(f_1^{u_1}(x), \ldots, f_m^{u_m}(x))\ \ \text{for all}\ \ x\in\R^n,$$ where $f_i^{u_i}(x):=f_i(x)-\langle u_i, x\rangle$, $i\in I$.  Consider the following perturbed 	optimization problem
\begin{equation}\label{problem-u}
	\mathrm{Min}\,_{\R^m_+}\{f^u(x)\;:\; x\in X\}, \tag{VP$_u$}
\end{equation}
where $u$ is the parameter of perturbation. The weak Pareto solution set of \eqref{problem-u} 
is denoted by Sol$^w(u)$. When $u=0$, one has $\mathrm{Sol}^w(0)=\mathrm{Sol}^w\eqref{problem}$.  

Solution stability is an interesting and very useful research field in optimization 
(see \cite{AT,Aus2,HL,rock-wet} among others) in virtue of its applications on 
concrete applications since, in practice, we are usually finding the solution of 
the optimization problem via numerical methods.

By a direct computation, one has  
\begin{align}
	& (f_i^{u_i})^{\infty} (y) = (f_i)^{\infty} (y) - \langle u_i, y \rangle, \label{equa-5a} \\
	& (f_i^{u_i})_q^{\infty} (y) = (f_i)_q^{\infty} (y) - \langle u_i, y \rangle. \label{equa-5q}
\end{align}

\begin{lemma}\label{CQ-u} If condition \eqref{equa-CQ} is satisfied, then there exists $\varepsilon>0$ such that for all $u=(u_1, \ldots, u_m)\in\R^{m\times n}$ satisfying $\|u\|:=\max_{i\in I}\|u_i\|<\varepsilon$, the following condition holds
	\begin{equation}\label{equa-6-new}
		X^{\infty} \cap\left(\bigcup_{i\in I}\mathcal{K}(f_i^{u_i})\right)=\{0\}.
	\end{equation}  
\end{lemma}	
\noindent{\it Proof}  
	It follows from  \eqref{equa-CQ} that $X^\infty\cap \mathcal{K}(f_i)=\{0\}$ for all $i\in I$. For each $i\in I$, we show that there exists $\varepsilon_i>0$ such that    
	\begin{equation}\label{equa-6}
		X^{\infty} \cap\mathcal{K}(f_i^{u_i})=\{0\} \ \ \forall u_i\in\mathbb{B}_{\varepsilon_i}.
	\end{equation}  
	Indeed, if otherwise, then for any  $k\in\N$, there is $u_i^k\in \mathbb{B}_{\frac{1}{k}}$ such that $$X^{\infty} \cap \mathcal{K}(f_i^{u^k_i})\neq \{0\}.$$ Hence, there exists $d_k\in X^{\infty} \setminus\{0\}$ such that $\big(f_i^{u^k_i}\big)^\infty (d_k)\leq 0$. Since $X^{\infty}$ is a closed cone, by passing a subsequence if necessary, we may assume that $h_k:=\frac{d_k}{\|d_k\|}$ converges to some $h\in X^{\infty}$ with $\|h\|=1$. For each $k>0$, by the positive homogeneity of $\big(f_i^{u^k_i}\big)^\infty$ and the fact that $\big(f_i^{u^k_i}\big)^\infty (d_k)\leq 0$, one has $\big(f_i^{u^k_i}\big)^\infty (h_k)\leq 0$.		By  \eqref{equa-5a}, we have
	\begin{equation*} 
		(f_i)^\infty (h_k) - \langle u_i^k, h_k \rangle= \big(f_i^{u^k_i}\big)^\infty (h_k)\leq 0.
	\end{equation*}
	Consequently,
	\begin{equation*} 
		(f_i)^\infty (h_k) \leq \langle u_i^k, h_k \rangle \ \ \forall k\in\N.
	\end{equation*}
	This and the lower semicontinuity of $(f_i)^\infty$  imply that
	\begin{equation*}
		(f_i)^{\infty} (h) \leq \liminf_{k \to \infty}(f_i)^{\infty} (h_k) \leq \lim_{k \to \infty} \langle u_i^k, h_k \rangle=0.  
	\end{equation*}
	Hence $h\in X^{\infty} \cap \mathcal{K}(f_i)$, which contradicts the fact that $X^{\infty} \cap \mathcal{K}(f_i)=\{0\}$. Thus \eqref{equa-6} holds. Put $\varepsilon:=\min\{\varepsilon_i\;:\; i\in I\}>0$. Then we have $X^{\infty} \cap\mathcal{K}(f_i^{u_i})=\{0\}$ for all $i\in I$ and $u_i\in\R^n$ satisfying   $\|u_i\|<\varepsilon$. Hence, \eqref{equa-6-new} is satisfied for all  $u=(u_1, \ldots, u_m)\in  \R^{m\times n}$ such that $\|u\|<\varepsilon$.
\qed
\begin{theorem}\label{theo:sta01} 
	Assume that condition \eqref{equa-CQ} is satisfied. 	 Then there exists $\varepsilon>0$ such that for all $u \in\R^{m\times n}$ satisfying $\|u\| <\varepsilon$, the following statements hold:
	\begin{enumerate} [$(a)$]
		\item $f^u$ is bounded from below on $X$.
		
		\item  For all $\lambda\in\Lambda$, the function $\varphi_\lambda^u(\cdot):=\langle \lambda, f^u(\cdot)\rangle$ is  coercive on $X$.
		
		\item  $\Sol^w(u)$ is nonempty and compact.
		
		\item  $\Limsup_{u \to 0}\Sol^w(u)\subset \Sol^w(0)$, i.e., the solution mapping $\Sol^w(\cdot)$ is closed at $0$.
		
		\item   $\Sol^w(\cdot)$ is usc at $0$.
	\end{enumerate}	
\end{theorem}
\noindent{\it Proof}  By Lemma \ref{CQ-u}, there exists $\varepsilon>0$ such that \eqref{equa-6-new} holds for all $u \in\R^{m\times n}$ satisfying $\|u\| <\varepsilon$.  Hence, conclusions $(a)$, $(b)$, and $(c)$ follow directly from Theorem \ref{weak-sharp-Thrm}.
	
	$(d)$: Let $\bar x$ be an arbitrary element belonging to $\Limsup_{u\to 0}\mathrm{Sol}^w\,(u)$. By definition,  there exist sequences $u^k:=(u_1^k, \ldots, u_m^k)\to 0$ and $x_k\in \mathrm{Sol}^w\,(u^k)$ such that $x_k\to \bar x$ as $k\to\infty$. Fix any  $x\in X$. For each $k\in \N$, since $x_k\in \mathrm{Sol}^w\,(u^k)$, there exists $i_k\in I$ such that
	\begin{equation*}
		f_{i_k}^{u^k_{i_k}}(x)-f_{i_k}^{u^k_{i_k}}(x_k)\geq 0.
	\end{equation*}
	Since $i_k\in I$ for all $k\in N$, there exists $i_0\in I$ such that $i_k=i_0$ for infinite many $k$. Hence, by passing to  subsequences if necessary, we may assume that
	\begin{equation*}
		f_{i_0}^{u^k_{i_0}}(x)-f_{i_0}^{u^k_{i_0}}(x_k)\geq 0 \ \ \forall k\in\N,
	\end{equation*}		
	or, equivalently,
	$$f_{i_0}(x_k) - \langle u_{i_0}^k, x_k \rangle\leq f_{i_0}(x) - \langle u_{i_0}^k, x \rangle \ \ \forall k\in\N.$$
	This and the lower semicontinuity of $f_{i_0}$ imply that
	\begin{align*}
		f_{i_0}(\bar x) \leq \liminf_{k \to \infty} f_{i_0}(x_k) & = \liminf_{k \to \infty} (f_{i_0}(x_k) - \langle u_{i_0}^k, x_k \rangle) \\
		& \leq \liminf_{k \to \infty} (f_{i_0}(x) - \langle u_{i_0}^k, x \rangle) = f_{i_0}(x).
	\end{align*}
	Hence, for each $x\in X$, there exists $i_0\in I$ such that $f_{i_0}(x)-f_{i_0}(\bar x)\geq 0$. This means that  $\bar x\in\Sol^w(0)$, as required.
	
	$(e)$:  Suppose on the contrary that  $\Sol^w(\cdot)$ is not usc at $0$. By definition, there exists an open set $U \subset \mathbb{R}^n$, with $\Sol^w(0) \subset U$, such that for every neighborhood $W$ of zero, there exists  $u\in W$ satisfying  $\Sol^w(u)\nsubseteq U$. Thus there exist   sequences $u^k:=(u_1^k, \ldots, u_m^k)\to 0$  and    $x_k \in \Sol^w(u^k) \setminus U$ for all $k \in \mathbb{N}$.  We consider two case of the sequence $x_k$ as follows.
	
	{\em Case 1.}  $\{x_k\}$ is bounded. Then, without loss of generality we may assume  that $x_k \to \bar{x}$.  By (d), one has  $\bar{x} \in \Sol^w(0) \subset U$.  However, since $x_k \notin U$ for all $k$ and $U$ is
	open, we see that $\bar x\notin U$, a contradiction. 
	
	{\em Case 2.} $\{x_k\}$ is unbounded. By passing to a subsequence if necessary, we may assume that $\|x_k\|\to\infty$ as $k\to\infty$. 	 Hence, without loss of generality one can assume that $\frac{x_k}{\|x_k\|}
	\to d \in X^{\infty}$ with $\|d\| = 1$. Fix any $x \in X$. Analysis similar to that in the proof of part $(d)$ shows that there exists $i_0\in I$ such that
	$$f_{i_0}(x_k) - \langle u_{i_0}^k, x_k \rangle\leq f_{i_0}(x) - \langle u_{i_0}^k, x \rangle \ \ \forall k\in\N.$$
We obtain that
	\begin{align*}
		(f_{i_0})^\infty (d)  \leq \liminf_{k \rightarrow \infty} \dfrac{f_{i_0}(x_k)}{\|x_k\|} &= \liminf_{k \rightarrow \infty} \dfrac{f_{i_0}(x_k) - \langle u_{i_0}^k, x_k \rangle}{\|x_k\|} 
		\\
		&\leq \liminf_{k \rightarrow \infty} \dfrac{f_{i_0}(x) - \langle u_{i_0}^k, x \rangle}{\|x_k\|} \, = 0.
	\end{align*} 
	Hence, $d\in X^{\infty} \cap \mathcal{K}(f_{i_0})$, contrary to \eqref{equa-CQ}. The proof is complete.
\qed

The following example  shows that there exists $u^k\to 0$ as $k\to \infty$ such that $\Sol^w(u^k)=\emptyset$ if condition \eqref{equa-CQ} is not satisfied.
\begin{example} \rm
	Consider the following vector optimization problem
	\begin{equation*}
		\mathrm{Min}_{\R^2_+}\,\{f(x):=(f_1(x),  f_2(x))  \;:\; x\in X\},
	\end{equation*}
	where $X=\R$, $f_1(x)=0$, and $f_2(x)=x$ for all $x\in\R$.		Then, $\Sol^w(0)=\R$ and 
	$$X^{\infty} \cap \big(\mathcal{K}(f_1)\cup \mathcal{K}(f_2)\big)=\R.$$
	Thus \eqref{equa-CQ} is not satisfied. For each $k$, let  $u^k=(-\frac{1}{k},\frac{1}{k})\in \mathbb{R}^2$. Then  $u^k\to 0$ as $k\to \infty$. Consider the following perturbed problem
	\begin{equation*}
		\mathrm{Min}_{\R^2_+}\,\bigg\{f^{u^k}(x):=\bigg(f_1^{-\frac{1}{k}}(x)=\frac{1}{k}x,  f_2^{\frac{1}{k}}(x)=x-\frac{1}{k}x\bigg)  \;:\; x\in X\bigg\}.
	\end{equation*}
	By definition, it is easy to see that $\Sol^w(u^k)=\emptyset.$
\end{example}

The following theorem presents sufficient conditions for the 
lower semicontinuity of the weak Pareto solution  map $\Sol^w(\cdot)$.

\begin{theorem}  \label{theorem12} 
	If the following conditions are satisfied:
	\begin{enumerate}[$(a)$]
		\item  $\Sol^w(0)$ is a singleton,
		
		\item  $X^{\infty} \cap \left(\bigcup_{i\in I}\mathcal{K}(f_i)\right)=\{0\}$, 
	\end{enumerate}
	then $\Sol^w(\cdot)$  is lsc  at $0$. 
\end{theorem}

\noindent{\it Proof}\ \ Suppose that $(a)$ and $(b)$ hold. By $(a)$, $\Sol^w(0) = \{\bar x\}$ for some
	$\bar{x} \in X$. Let $U$ be an open neighborhood containing  
	$\bar x$. Then, by (b) and Lemma \ref{CQ-u}, there exists $\varepsilon_1>0$ such that for all $u=(u_1, \ldots, u_m)\in\R^{m\times n}$ satisfying $\|u\|<\varepsilon_1$, condition  \eqref{equa-6-new} holds.
	It follows from Theorem \ref{theo:sta01}(c) that  $\Sol^w(u) \neq \emptyset$ for every $u$ with $\Vert u \Vert < \varepsilon_1$. Since $\Sol^w(\cdot) $ is usc at $0$, by Theorem \ref{theo:sta01}(e), there exists $\varepsilon_{2} > 0$ such that $\Sol^w(u) \subset U$ for every $u$ satisfying $\Vert u \Vert < \varepsilon_{2}$.  
	By taking $\varepsilon := \min\{\varepsilon_{1}, \varepsilon_{2}\} > 0$, we obtain that 	$\Sol^w(u) \cap U \neq \emptyset$ for all $u$ satisfying $\|u\|<\varepsilon$. Hence  $\Sol^w(\cdot)$ is lsc at $0$. 	
\qed

\begin{theorem}\label{t1}
	Assume that  $X$ is convex and $f$ is $\R^m_+$-convex. If $\Sol^w(\cdot)$ is lsc at $0$, then the following statements are equivalent:
	\begin{enumerate} [$(a)$]
		\item $\Sol^w(0)$ is bounded.
		\item $X^{\infty} \cap \left(\bigcup_{i\in I}\mathcal{K}(f_i)\right)=\{0\}.$
	\end{enumerate}
\end{theorem}
\noindent{\it Proof} 
	By assumptions, we first claim that
	\begin{equation*}
		X^{\infty} \cap \left(\bigcap_{i\in I}\mathcal{K}(f_i)\right)=\{0\}.
	\end{equation*}
	Indeed, if otherwise, then there exists a nonzero vector $d \in X^{\infty}$ such that $(f_i)^{\infty} (d) \leq 0$ for all $i\in I$. For each $k\in\mathbb{N}$, let $u_i^k:=\frac{1}{k}d$ for all $i\in I$. Then, $u_i^k \rightarrow 0$ as $k\to\infty$ and $\langle u_i^k, d \rangle > 0$ for all $k \in \mathbb{N}$ for all $i\in I$. Hence,  
	$$\big(f_i^{u_i^k}\big)^\infty(d)= (f_i)^{\infty} (d) - \langle u_i^k, d 
	\rangle < 0 \ \  \forall k \in \mathbb{N}.$$
	We show that $\Sol^w(u^k)$ is empty for every $k\in\N$ and so which contradicts the lower semicontinuity of $\Sol^w(\cdot)$ at $0$ and the fact that $u^k\to 0$ as $k\to\infty$.  Indeed, if otherwise, then there exists $k_0\in\N$ such that $\Sol^w(u^{k_0}) \neq \emptyset$. Fix any $x_0\in \Sol^w(u^{k_0})$. It follows from the convexity of $X$ that $x_0+td\in X$ for all $t>0$.  By \cite[Proposition 2.5.2]{AT} and the convexity of  $f_i$, one has
	$$\big(f_i^{u_i^{k_0}}\big)^\infty(d)=\lim_{t\to\infty}\dfrac{f_i^{u_i^{k_0}}(x_0+td)-f_i^{u_i^{k_0}}(x_0)}{t}= \lim_{t\to\infty}\dfrac{f_i^{u_i^{k_0}}(x_0+td)}{t}.$$
	Hence,
	$$\lim_{t\to\infty}\dfrac{f_i^{u_i^{k_0}}(x_0+td)}{t}=\big(f_i^{u_i^{k_0}}\big)^\infty(d)<0.$$
	This implies that $\lim_{t\to\infty} {f_i^{u_i^{k_0}}(x_0+td)}=-\infty$ for all $i\in I$. Hence, there exists $t$ large enough such that 
	\begin{equation*}
		f^{u^{k_0}}(x_0+td)\in f^{u^{k_0}}(x_0)-\mathrm{int}\,\R^m_+,
	\end{equation*}
	contrary to the fact that $x_0\in \Sol^w(u^{k_0})$, as required.
	
	Now, by the lower semicontinuity of $\Sol^w(\cdot)$ at $0$, we see that $\Sol^w(0)$ is nonempty. Hence the statement $(a)$  is equivalent to  the nonemptiness and compactness of $\Sol^w(0)$. The equivalence of $(a)$ and $(b)$ follows immediately from   \cite[Theorem 3.1]{Wang-Fu-24}. The proof is complete.
\qed

The following example not only illustrates Theorem \ref{t1} but also shows that the lower semicontinuity of $\Sol^w(\cdot)$ does not imply the singletoness of the set  $\Sol^w(0)$.
\begin{example} \rm 
	Consider the following vector optimization problem
	\begin{equation*}
		\mathrm{Min}_{\R^2_+}\,\{f(x):=(f_1(x),  f_2(x))  :\ x\in X\},
	\end{equation*}
	where $X=[-1,+\infty)$, $f_1(x)=x$, and $f_2(x)=x^2$ for all $x\in\R$. Then, $\Sol^w(0)=[-1,0]$ is not a singleton and 
	$$X^{\infty} \cap \big(\mathcal{K}(f_1)\cup \mathcal{K}(f_2)\big)=\{0\}.$$
	For each $u=(u_1,u_2)\in \mathbb{R}^2$, we have   the following perturbed problem
	\begin{equation*}
		\mathrm{Min}_{\R^2_+}\,\{f^u(x):=(f_1^{u_1}(x)=x-u_1x,  f_2^{u_2}(x)=x^2-u_2x)\;:\; x\in X=[-1,+\infty)\}.
	\end{equation*}
	We can check that $$\Sol^w(u)=\bigg[-1,\frac{u_2}{2(1-u_1)}\bigg]$$
	and $\Sol^w(\cdot)$ is lsc at $0$.
\end{example}

\section{The Quasiconvex Case}\label{section5}

In this section, we apply our previous results to the particular case when the 
objective function in problem \eqref{problem} is quasiconvex or $\alpha$-robustly
quasiconvex (see Definition \ref{alpha:robust}).

As mentioned in the introduction, when the function $f$ is nonconvex, then the 
usual asymptotic function $f^{\infty}$ does not provide adequate information
on the behavior of $f$ at infinity. For instance, and in relation to Theorem 
\ref{weak-sharp-Thrm}, we mention that when the function $f$ is quasiconvex, 
the assumption 
$$X^{\infty} \cap \mathcal{K}(f) =  \{0\},$$
is too restrictive. Indeed, let us consider the one-dimensional real-valued 
function $f: \mathbb{R} \rightarrow \mathbb{R}$ given by $f(x) = \sqrt{\lvert 
	x \rvert}$ and $X=\mathbb{R}$. Here $f$ is coercive and ${\rm argmin}_{X}\,f 
= \{0\}$ is a singleton. However, $X^{\infty} \cap \mathcal{K}(f) = \mathbb{R}$ 
and hence, Theorem \ref{weak-sharp-Thrm} cannot be applied  to this basic situation.

On the other hand, if we use any of the generalized asymptotic functions
$f^{\infty}_{q}$ and $f^{\infty}_{\lambda}$, we obtain that
$$f^{\infty}_{q} (u) = + \infty\ \  \forall  u \neq 0,$$ and for every $\lambda\in\R$ such that $\text{lev}\,(f,\lambda)\neq\emptyset$, 
$$0 < f^{\infty}_{\lambda} (u) \leq 2\ \  \forall   u \neq 0,$$
where
$$f^{\infty}_{\lambda} (u):=\sup_{x\in\text{lev}\,(f,\lambda)}\;\sup_{t>0}\frac{f(x+tu)-\lambda}{t}.$$
Therefore, 
$$X^{\infty} \cap  \mathcal{K}_{q} (f) =  \{0\}, ~~ {\rm and}~~ X^{\infty} \cap 
\mathcal{K}_{\lambda} (f) =  \{0\},$$
where $\mathcal{K}_{q}(f) := \{d \in \mathbb{R}^{n}: \, f^{\infty}_{q} (d) \leq 0\}$ and
$\mathcal{K}_{\lambda}(f) := \{d \in \mathbb{R}^{n}: \, f^{\infty}_{\lambda} (d) \leq 0\}$,
respectively. 

Before stating main results of this section, we give a result on the nonemptiness 
and compactness of the solution set to constrained optimization problems by 
using the $q$-asymptotic function.

\begin{lemma}\label{lemma-51} 
	Assume that $X$ is  convex  and   $f$  is   $\mathbb{R}^m_+$-quasiconvex. Then, the  following assertions are equivalent:
	\begin{enumerate}[$(a)$]
		\item $\Sol^w\eqref{problem}$ is nonempty and compact.
		
		\item $X^\infty\cap \left(\bigcup_{i\in I}\mathcal{K}_q(f_i)\right)=\{0\}$.
	\end{enumerate}
\end{lemma} 

\noindent{\it Proof} \ \ 
	In order to prove the equivalence, we first show that
	\begin{equation*}
		(f_i+\delta_X)_q^\infty(d)=(f_i)_q^\infty(d)+\delta_{X^\infty}(d)\ \ \forall i\in I.
	\end{equation*}
	Indeed, by definition of the $q$-asymptotic function and the convexity of $X$, we have
	\begin{align*}
		(f_i+\delta_X)_q^\infty(d)&=\sup_{x \in {\rm dom}\,(f_i+\delta_X)} \sup_{t>0} \frac{(f_i+\delta_X)(x+td) - (f_i+\delta_X)(x)}{t}
		\\
		&=\sup_{x \in X} \sup_{t>0} \frac{f_i(x+td)+\delta_X(x+td) - f_i(x)}{t}
		\\
		&=
		\begin{cases}
			(f_i)_q^\infty(d), &\text{if}\ \ d\in X^\infty,
			\\
			+\infty, &\text{otherwise},
		\end{cases}
		\\
		&=(f_i)_q^\infty(d)+\delta_{X^\infty}(d).
	\end{align*}		
	Hence, 
	\begin{align*}
		\mathcal{K}_q(f_i+\delta_X)=\{d\in\mathbb{R}^n  \;:\;(f_i)_q^\infty(d)+\delta_{X^\infty}(d)\leq 0\}=X^\infty\cap \mathcal{K}_q(f_i) \ \ \forall i\in I.
	\end{align*}		
	This implies that
	\begin{equation*}
		\bigcup_{i\in I}\mathcal{K}_q(f_i+\delta_X)=X^\infty\cap\left(\bigcup_{i\in I}\mathcal{K}_q(f_i)\right).
	\end{equation*}		
	Combining this with  \cite[Proposition 4.5]{F-L-Vera} we obtain the conclusion of the lemma. 
\qed

In the next proposition we improve Theorem \ref{weak-sharp-Thrm} for proper, 
lsc, $\mathbb{R}^m_+$-$\alpha$-robustly quasiconvex functions. 

\begin{proposition}
	Let $X$ be a convex set and $f$ be a proper, 
	lsc, $\mathbb{R}^m_+$-$\alpha$-robustly quasiconvex mapping
	$(\alpha>0)$. If 
	\begin{equation}\label{equa-CQq}
		X^\infty\cap\left(\bigcup_{i\in I}\mathcal{K}_q(f_i)\right) = \{0\},
	\end{equation}
	then there exists $\varepsilon > 0$ such that for all $u \in \mathbb{R}^{m\times n}$ with $\|u\|<\varepsilon$,
	the following statements hold:
	\begin{enumerate} [$(a)$]
		\item $f_u$ is bounded from below on $X$.
		
		\item ${\rm \Sol}^w(u)$ is nonempty and compact.
		
		\item $ \Limsup_{u \rightarrow 0} {\rm \Sol}^w (u) \subset {\rm \Sol}^w (0)$, i.e., the solution mapping $\Sol^w(\cdot)$ is closed at $0$.
		
		\item $\Sol^w(\cdot)$ is usc at $0$.
	\end{enumerate}
\end{proposition}

\noindent{\it Proof} 
	We first show that there exists $\varepsilon > 0$ such that  the following condition holds  
	\begin{equation}\label{q:assump}
		X^{\infty} \cap \left(\bigcup_{i\in I}\mathcal{K}_q(f^{u_i}_i)\right)= \{0\} \ \ \forall u \in 
		\mathbb{R}^{m\times n}, \|u\|<\varepsilon.
	\end{equation}
	Indeed, suppose on the contrary that for every $k \in \mathbb{N}$, there 
	exists $u^{k} \in \mathbb{R}^{m\times n}$ with $\|u^k\|<\frac{1}{k}$ such that condition \eqref{q:assump} does not hold, i.e., there exist $i_k\in I$ and $d_{k} \in X^{\infty} 
	\backslash \{0\}$ such that 
	$$\bigg(f_{i_{k}}^{u_{i_k}^k}\bigg)^{\infty}_{q} (d_{k}) \leq 0.$$
	By the finiteness of $I$, there exists $i_0\in I$ such that the inequality
	\begin{equation}\label{equa-new-1}
		\bigg(f_{i_{0}}^{u_{i_0}^k}\bigg)^{\infty}_{q} (d_{k}) \leq 0
	\end{equation}
	holds true for infinite many $k\in\N$. By passing to subsequences if necessary we may assume that \eqref{equa-new-1} is satisfied for all $k\in\mathbb{N}$. 
	
	Since $X^{\infty}$ is a closed cone, by passing a subsequence if necessary, we 
	may assume that $\left\{h_{k} := \frac{d_{k}}{\lVert d_{k} \rVert}\right\}_{k}
	\subset X^{\infty}$ converges to $h \in X^{\infty}$ with $\lVert h \rVert = 1$.
	
	For every $k \in \mathbb{N}$, since $\bigg(f_{i_{0}}^{u_{i_0}^k}\bigg)^{\infty}_{q}$ is positively homogeneous of degree one and $\bigg(f_{i_{0}}^{u_{i_0}^k}\bigg)^{\infty}_{q}(h_{k}) \leq 0$, and by using relation \eqref{equa-5q}, we have
	\begin{align*}
		\bigg(f_{i_{0}}^{u_{i_0}^k}\bigg)^{\infty}_{q}(h_{k})	= (f_{i_0})^{\infty}_{q} (h_{k}) - \langle u^{k},  h_{k} \rangle\leq 0. 
	\end{align*}
	Since $f_{i_0}$ is lsc, $(f_{i_0})^{\infty}_{q}$ is lsc too by \cite[p. 118]{FFB-Vera}, thus
	$$(f_{i_0})^{\infty}_{q} (h) \leq \liminf_{k \rightarrow + \infty} (f_{i_0})^{\infty}_{q} (h_{k}) \leq 
	\liminf_{k \rightarrow + \infty} \, \langle u^{k},  h_{k} \rangle = 0.$$
	Hence, $h \in X^{\infty} \cap \mathcal{K}_{q} (f_{i_0})$ with $h \neq 0$, a contradiction. 
	Therefore, relation \eqref{q:assump} holds.
	
	Since $X$ is convex and $f$ is $\R^m_+$-$\alpha$-robustly quasiconvex, we obtain that the function $f_{u}(x)=f(x) - \langle u, x \rangle$ is $\R^m_+$-quasiconvex for all $u\in \mathbb{R}^{m\times n}$ with $\|u\|<\alpha$. Let $\varepsilon$ be satisfied \eqref{q:assump} and $\varepsilon<\alpha$. Hence, $(b)$ follows from relation
	\eqref{q:assump} and Lemma \ref{lemma-51}. 
	
	The relation \eqref{q:assump} gives
	\begin{equation*}
		X^{\infty} \cap \mathcal{K}_q(f^{u_i}_i)= \{0\} \ \ \forall u \in 
		\mathbb{R}^{m\times n}, \|u\|<\varepsilon.
	\end{equation*}
	Combining this with \cite[Proposition 5.4]{Lara-Tuyen-Nghi} yields that $f^{u_i}_i$ is bounded from below on $X$ and $(a)$ follows.
	
	The proof of part $(c)$ is quite similar to the proof of Theorem \ref{theo:sta01}(d), so omitted.
	
	$(d)$: Suppose on the contrary that  $\Sol^w(\cdot)$ is not usc at $0$. Then, there exist an open set $U \subset \mathbb{R}^n$, with $\Sol^w(0) \subset U$,  sequences $u^k$ and $x_k$ such that  $u^k\rightarrow 0$,  $x_k \in \Sol^w(u^k)\setminus U$ for all $k \in \mathbb{N}$.   We consider two case of the sequence $x_k$ as follows.
	
	{\em Case 1.}  $\{x_k\}$ is bounded. Then, without loss of generality we may assume  that $x_k \to \bar{x}$.  By $(c)$, one has  $\bar{x} \in \Sol^w(0) \subset U$.  However, since $x_k \notin U$ for all $k$ and $U$ is
	open, we see that $\bar x\notin U$, a contradiction. 
	
	{\em Case 2.} $\{x_k\}$ is unbounded. By passing to a subsequence if necessary, we may assume that $\|x_k\|\to\infty$ as $k\to\infty$. 	 Hence, without loss of generality one can assume that $\frac{x_k}{\|x_k\|}
	\to d \in X^{\infty}$ with $\|d\| = 1$. Fix any $y \in X$. Analysis similar to that in the proof of Theorem \ref{theo:sta01}(d)   shows that  there exists $i_0\in I$ such that
	\begin{equation}\label{e1}
		f_{i_0}(x_k) - \langle u_{i_0}^k, x_k \rangle\leq f_{i_0}(y) - \langle u_{i_0}^k, y \rangle \ \ \forall k\in\N.
	\end{equation}
	Since $u^{k} \to 0$, there exists $k_{1} \in \mathbb{N}$ such that 
	$\|u^{k}\|<\alpha$ for all $k \geq k_{1}$. Furthermore, since $\|x_k\| \to + \infty$  for every $t>0$, there exists $k_{2} \in \mathbb{N}$ such that $0 < \frac{t}{\|x_k\|} < 1$ for all $k \geq k_{2}$. 
	
	Since $X$ is convex and $f$ is $\mathbb{R}^m_+$-$\alpha$-robustly quasiconvex, by \eqref{e1}, we obtain for every $k \geq k_{0} := \max\{k_{1}, k_{2}\}$ that
	\begin{align*}
		f_{i_0}\left( \left(1 - \frac{t}{\|x_k\|}\right) y + \frac{t}{\|x_k\|} x_{k} \right) -&\left\langle u^{k}_{i_0}, \left(1 - \frac{t}{\|x_k\|}\right) y + \frac{t}{\|x_k\|} x_{k}  \right\rangle\\
			=&f^{u^k_{i_0}}_{i_0}\left( \left(1 - \frac{t}{\|x_k\|}\right) y + \frac{t}{\|x_k\|} x_{k} \right)  
		\\ 
		\leq&\max\{f_{i_0}^{u^{k}_{i_0}} (x_{k}),  f_{i_0}^{u^{k}_{i_0}} (y)\}\\
		=& f_{i_0}(y) - \langle u^{k}_{i_0}, y \rangle.
	\end{align*}
	Letting $k\to \infty$, we have
	$$f_{i_0}(y + td) \leq f_{i_0}(y)  \ \  \forall ~ t > 0,\ \  \forall  y \in X,$$
	due to the lower semicontinuity  of  $f_{i_0}$.
	This implies that 
	$(f_{i_0})^{\infty}_{q} (d) \leq 0$ and we therefore get $d \in  
	X^\infty\cap \mathcal{K}_q(f_{i_0})$, contrary to \eqref{equa-CQq}. The proof is complete.
\qed

Furthermore, we also adapt the results for solution stability for the quasiconvex 
case below.

\begin{proposition} 
	Let $X$ be a convex set and $f$ be an  $\mathbb{R}^m_+$-$\alpha$-robustly quasiconvex mapping ($\alpha 
	> 0$). If the following conditions hold:
	\begin{itemize}
		\item[$(a)$] $\Sol^w(0)$ is a singleton;
		
		\item[$(b)$] $X^\infty\cap\left(\bigcup_{i\in I}\mathcal{K}_q(f_i)\right)= \{0\}$;
	\end{itemize}
	then $\Sol^w(\cdot)$ is lsc at $0$. Conversely, if  $\Sol^w(\cdot)$ is lsc at $0$, then  the following condition is satisfied
	\begin{equation}\label{equa-new-2}
		X^\infty\cap\left(\bigcap_{i\in I}\mathcal{K}_q(f_i)\right)= \{0\}.
	\end{equation}
\end{proposition}

\noindent{\it Proof} 
	The proof of the first conclusion is quite similar to  that of the proof of Theorem \ref{theorem12}, so omitted. We now suppose on the contrary that $\Sol^w(\cdot)$ is lsc at $0$ but \eqref{equa-new-2} does not hold. Then, $\Sol^w(0)$ is nonempty and there exists $d\in X^\infty\setminus\{0\}$ such that
	$(f_{i})^\infty_q(d)\leq 0$  for every  $i\in I$. Let  $\bar x\in\Sol^w(0)$. Clearly, $\bar x+td\in X$ for all $t>0$. Let $u^k: =(\frac{1}{k}d, \ldots, \frac{1}{k}d) $. Then $u^k\to 0$ as $k\to\infty$ and for each $i\in I$,
	$$(f_{i}^{u^k_i})^\infty_q(d)=(f_{i})^\infty_q(d)-\langle u^k_i, d\rangle=(f_{i})^\infty_q(d)-\frac{1}{k}\|d\|^2<0.$$
	By definition, we have
	\begin{equation*}
		\sup_{t>0}\frac{f_{i}^{u^k_i}(\bar x+td)-f_{i}^{u^k_i}(\bar x)}{t}<0.   
	\end{equation*}
	This implies that $f_{i}^{u^k_i}(\bar x+td)\to-\infty$ as $t\to\infty$ for every $i\in I$. This means that for each $z\in X$ there exists $t>0$ large enough such that
	$$\max_{i\in I}\{f_{i}^{u^k_i}(\bar x+td)-f_{i}^{u^k_i}(z)\}<0.$$
	It follows that  $\Sol^w(u^k)=\emptyset$. Hence, $\Sol^w(\cdot)$ is not lsc at $0$, a contradiction. The proof is complete.
\qed

\section{Conclusions}

In this paper,  we have used the tool of asymptotic analysis to investigate  existence and stability of weak Pareto solutions of constrained  nonconvex vector optimization problems.  Under the proposed assumption that 
$$X^{\infty} \cap \bigg(\bigcup_{i\in I}\mathcal{K}(f_i)\bigg)=\{0\},$$ 
which is a relation between the asymptotic of  objective functions and the asymptotic cone of the constraint set,  we have presented new results on existence of weak Pareto solutions, the weak sharp minima at infinity property,  the upper/lower semicontinuity, the nonemptiness and compactness, and the closedness  of  the weak Pareto  solution map of the considered problem with linear perturbation. We have used $q$-asymptotic functions to obtain  similar results to the particular case when the 
objective function   is quasiconvex or $\alpha$-robustly
quasiconvex.
 
\section*{Acknowledgments} We are grateful to two anonymous referees  for their careful readings and valuable suggestions which improved the presentation of this manuscript. In particular, we would like to thank Referee $\#$2 for his/her suggestions, which helped us to improve the proof of Theorem \ref{weak-sharp-Thrm}(e). 

\section*{Funding} 
This research is funded by the Hanoi Pedagogical University 2 [grant number HPU2.2025-UT-04].

\section*{Disclosure statement} 
The authors declare that they have no conflict of interest.

\section*{Data availability} There is no data included in this paper.

\bibliographystyle{amsplain}

\end{document}